\documentclass[a4paper]{jpconf}
\usepackage{graphicx}

\newenvironment{dedication}
  { \itshape             % the text is in italics
 \raggedleft
  }
  {\par % end the paragraph
  }

\usepackage{amsmath,amssymb,cite,comment}

\def\md{\medskip}
\def\nn{\nonumber}
\def\bpm{\begin{pmatrix}}
\def\epm{\end{pmatrix}}
\def\D{{\Delta}}\def\veps{\varepsilon}
\def\cgc{{\cg^\bbc}} \def\r{{\rho}}
\def\bca{\begin{cases}}
\def\eca{\end{cases}}

\def\rf#1#2{(\ref{#1}{#2})}

\def\ha{{\textstyle{1\over2}}}
\def\trha{{\textstyle{3\over2}}}
\def\third{{\textstyle{1\over3}}}

\def\nt{\noindent}
\input epsf.tex
\newcount\figno
\def\fig#1#2#3{
\par\begingroup\parindent=0pt\leftskip=1cm\rightskip=1cm\parindent=0pt
\baselineskip=11pt \global\advance\figno by 1 %\midinsert
\epsfxsize=#3 \centerline{\epsfbox{#2}} \vskip 12pt #1\par
\endgroup\par}
\def\figlabel#1{\xdef#1{\the\figno}}
\def\L{\Lambda} \def\l{\lambda} 
\def\ca{{\cal A}} \def\cb{{\cal B}} 
  \def\cf{{\cal F}}
\def\cg{{\cal G}} \def\ch{{\cal H}} 
 \def\ck{{\cal K}} \def\cl{{\cal L}}
\def\cm{{\cal M}} \def\cn{{\cal N}} 
\def\cp{{\cal P}} \def\cq{{\cal Q}}  \def\cs{{\cal S}}

\def\a{\alpha} 
\def\b{\beta} 
\def\d{\delta}

\def\bbz{\mathbb{Z}}%{Z\!\!\!Z}
\def\bbc{\mathbb{C}}%{I\!\!\!\!C}}
\def\bbr{\mathbb{R}}%{{I\!\!R}}
\def\bbn{\mathbb{N}}%{I\!\!N}

\def\eqn#1{\begin{equation}\label{#1}}
\def\ee{\end{equation}}

\def\bea{\begin{eqnarray}}
\def\eea{\end{eqnarray}}

\def\eqnn#1{\begin{eqnarray}\label{#1}}

\newcommand{\eqna}[1]{\begin{subequations} \label{#1}
\begin{eqnarray}}
\def\eena{\end{eqnarray}
\end{subequations}}

 \usepackage{graphicx}

\begin{document}

\title{Multiplet Classification of Reducible Verma Modules over the $G_2$ Algebra}

% \title{Invariant Differential Operators   of the $G_2$ Algebra:\\
% Classification of Verma Modules}

\author{V.K. Dobrev}

\address{Institute for Nuclear
Research and Nuclear Energy,
 Bulgarian Academy of Sciences,  72
Tsarigradsko Chaussee, 1784 Sofia, Bulgaria; dobrev@inrne.bas.bg}

%\ead{dobrev@inrne.bas.bg}
%{\it dobrev@inrne.bas.bg}

 \begin{abstract}
In the present paper we continue the project of systematic
construction of invariant differential operators on the example of
the non-compact  algebra  $G_{2(2)}$ which is split real form of $G_2$.  We give
 the classification of reducible Verma modules $G_2$. We give also the singular vectors between
 these modules, thus setting the stage for construction of the  invariant differential operators over
 $G_{2(2)}$.
\end{abstract}

\begin{dedication}
Dedicated to I.E. Segal (1918-1998)  in commemoration of the centenary of his birth.\\
The author remembers with great pleasure the talk he gave at Segal's seminar at MIT in 1975.
\end{dedication}

\md

\section{Introduction}

Invariant differential operators   play very important role in the
description of physical symmetries. The general scheme for constructing these
operators was given some time ago \cite{Dob}.
In recent papers \cite{Dobinv,Dobparab} we started the systematic explicit
construction of the invariant differential operators.

The first task in the construction is to make the multiplet classification of
the reducible Verma modules over the algebra in consideration following \cite{Dobmul}.
Such classification provides the weights of embeddings between the Verma modules via
the singular vectors, and thus, by \cite{Dob}, the weights of the invariant differential operators.

We have done the multiplet classification for many real non-compact algebras, first from the class
of algebras that have discrete series representations, see \cite{Dobk1}.
In the present paper we  focus on the complex algebra ~$G_2$~ and on its
 split real form algebra  $G_{2(2)}$. We present the multiplet classification of
 the reducible of  Verma modules over $G_2\,$. We give also the singular vectors between
 these modules. By the scheme of \cite{Dob} these explicit expressions produce
 the invariant differential operators.

 This paper is a   sequel of \cite{Dobinv} and \cite{Dobparab}
  and due to the lack of space we refer to these papers
 for motivations and \cite{Dobk1} for extensive list of literature on the subject.
 For other approaches and applications of $G_2$, see, e.g., \cite{G2refs}.

\section{Preliminaries}

\subsection{Lie algebra}

We start with the complex Lie algebra  ~$\cg^\bbc ~=~ G_2$. We use
the standard definition of ~$\cg^\bbc$~ given in terms of the
Chevalley generators $X^\pm_i ~, ~H_i ~, ~i=1,2$, by the relations~:
\eqnn{com}
 &[H_j\,, \,H_k] \, = \, 0 \,, \,\,\,[H_j\,, \,X^\pm_k] \, = \, \pm
a_{jk} X^\pm_k \,,
\,\,\,[X^+_j \,, \,X^-_k] \, = \,
\d_{jk} \,H_j \,,  \\
&\sum_{m=0}^n \,(-1)^m \,\left({n \atop m}\right)
\,\left(X^\pm_j\right)^m \,X^\pm_k \,\left(X^\pm_j\right)^{n-m}
\,=\, 0 \,, \,\,j \neq k \,, \,\,n = 1 - a_{jk} \,, \nn\eea where
\eqn{cart}  (a_{jk}) \,=\, (\a^\vee_j , \a_k) \,=\,
\bpm 2 & -1\cr -3 &2 \epm\ee is the Cartan matrix of $\cg^\bbc$, \,\,$\a^\vee_j
\,\equiv\, {2 \a_j\over (\a_j , \a_j)}$\, is the co-root of
$\a_j\,$, \,\, $(\cdot , \cdot)$ is the scalar product of the roots,
so that the nonzero products between the simple roots are:
$(\a_1,\a_1) = 6$, $(\a_2,\a_2) = 2$, $(\a_1,\a_2) = -3$. The
elements $H_i$ span the Cartan subalgebra $\ch$ of $\cg^\bbc$, while
the elements $X^\pm_i$ generate the subalgebras $\cg^\pm$. We shall
use the standard triangular decomposition \eqn{deca} \cg^\bbc \,=\,
\cg_+\otimes \ch\oplus \cg_- \,, \,\,\cg_\pm \,\equiv
\,\mathop{\oplus}\limits_{\a\in\D^\pm} \,\cg_\a \,, \ee where
$\D^+$, $\D^-$, are the sets of positive, negative, roots, resp.
Explicitly    we have:
\eqn{posg2} \D^+ = \{ \a_1, ~\a_2, ~\a_1 + \a_2, ~\a_1 +
2\a_2, ~\a_1 + 3\a_2, ~2\a_1 + 3\a_2 \} \ee

Thus, $G_2$ is 14--dimensional ($14 = \vert \D\vert +$ rank $G_2$).

For the simple roots we may choose in terms of ortho-normal basis $\veps_1,\veps_2,\veps_3$~:
\eqn{pig2} \a_1 = \veps_1 +
\veps_3 - 2\veps_2 , \ \ \ \a_2 = \veps_2 - \veps_3 \ee

For future reference we introduce notation for the non-simple roots:
\eqnn{nons} &&\a_3 \equiv \a_1 + \a_2 = \veps_1 - \veps_2,\quad  \a_4 \equiv \a_1 +
2\a_2 = \veps_1 - \veps_3, \\
&&\a_5 \equiv \a_1 + 3\a_2 = \veps_1 +
\veps_2 - 2\veps_3, \quad \a_6\equiv 2\a_1 + 3\a_2 = 2\veps_1 -
\veps_2 - \veps_3 \nn\eea

With the chosen normalization the roots $\a_1,~\a_5=\a_1 + 3\a_2,~ \a_6=2\a_1 +
3\a_2$ have length 6, while $\a_2,~\a_3=\a_1 + \a_2,~\a_4=\a_1 + 2\a_2$ have
length 2. The dual roots are:
\eqnn{dual} && \a_1^\vee ~=~ \a_1/3, \qquad \a_2^\vee
~=~ \a_2, \nn \\ && \a_3^\vee=(\a_1 + \a_2)^\vee ~=~ \a_1 + \a_2 ~=~ 3\a_1^\vee +
\a_2^\vee, \nn \\ && \a_4^\vee=(\a_1 + 2\a_2)^\vee ~=~ \a_1 + 2\a_2 ~=~ 3\a_1^\vee +
2\a_2^\vee, \nn \\ && \a_5^\vee=(\a_1 + 3\a_2)^\vee ~=~ (\a_1 + 3\a_2)/3 ~=~
\a_1^\vee + \a_2^\vee, \nn \\ &&\a_6^\vee=(2\a_1 + 3\a_2)^\vee ~=~ (2\a_1 + 3\a_2)/3
~=~ 2\a_1^\vee + \a_2^\vee \eea

Note that the roots $\a_2,\a_3,\a_4$ form the $sl(3)$ root system, the first two being
  the two simple roots. The roots $\a_1,\a_5,\a_6$ also form the $sl(3)$ root
  system, the standard normalization being achieved after rescaling each root by factor $\sqrt{3}$.
  Note also the cases: ~$(\a_1,\a_4)=0$, ~$(\a_2,\a_6)=0$, ~$(\a_3,\a_5)=0$.

The Weyl group ~$W(G_2)$~ of ~$G_2$~ is the dihedral group of order 12. This
follows from the fact that ~$(s_1\,s_2)^6=1$, where ~$s_1\,,s_2$~
are the two simple reflections. Using the general formula:
\eqnn{sapr}  s_\a (\lambda) ~=~ \lambda - {2
\langle \lambda, \a{\rangle}\over \langle\a,\a{\rangle}} \a ~=~
\lambda - \langle \lambda, \a^\vee{\rangle} \a ~,  \nn\\  \lambda
\in \ch^*, ~\a \in \Delta , ~~~\a^\vee \equiv {2 \over
\langle\a,\a{\rangle}} \a ~, \eea
we have the following action of the simple reflections: ~$s_i \equiv s_{\a_i}$~:
\eqnn{scat}
&& s_1 (\a_1,\a_2,\a_3,\a_4,\a_5,\a_6) =
 (-\a_1,\a_3,\a_2,\a_4,\a_6,\a_5) \\
&& s_2 (\a_1,\a_2,\a_3,\a_4,\a_5,\a_6) =
 (\a_5,-\a_2, \a_4, \a_3, \a_1,\a_6)
\eea
The 12 elements of ~$W(G_2)$~ are given in terms of the simple reflections as follows:
\eqnn{weylg} && W(G_2) = \{ 1,~ s_1,~ s_2,~ s_1 s_2,~  s_2 s_1,~
s_1 s_2s_1,~  s_2 s_1 s_2 ,~ s_1 s_2s_1 s_2 ,~ s_2 s_1s_2 s_1, \nn\\
&& s_1 s_2s_1 s_2 s_1,~ s_2 s_1s_2 s_1 s_2,~
~ s_1 s_2 s_1s_2 s_1 s_2 = ~ s_2 s_1s_2 s_1 s_2 s_1\}
\eea
Note the expressions for reflections corresponding to non-simple roots:
\eqn{reflns} s_{\a_3} ~=~ s_1 s_2 s_1 \ , \quad s_{\a_4} ~=~ s_2 s_1 s_2 s_1 s_2  \ , \quad
 s_{\a_5} ~=~ s_2 s_1 s_2 \ , \quad s_{\a_6} ~=~  s_1 s_2 s_1 s_2 s_1 \ee

Let us denote the root space vector of \,$\cg_\a$\, by \,$X_\a$, or more explicitly:
\,$X^\pm_k \,\equiv\, X_{\pm\a_k}\,$, \,$k=1,...,6$. To give the full
Cartan-Weyl basis we need to define also $X_k^\pm\,$, \,$k=3,...,6$, for which we follow \cite{Dob}:
\eqnn{cww} && X^\pm_3 = \pm [X^\pm_1, X^\pm_2], \qquad X^\pm_4 = \pm \, [  X^\pm_2, X^\pm_3]
 \ , \\ && X^\pm_5 = \pm \,\frac{1}{\sqrt{3}} [  X^\pm_4, X^\pm_2] \ , \qquad X^\pm_6
 = \pm \, [  X^\pm_1, X^\pm_5]
   \nn\eea
  Then we have for ~$H_k \,\equiv \,[ X^+_k , X^-_k ]$:
\eqn{crr}   H_3  \,=\, 3H_1 + H_2 \,, \quad
H_4  \,=\, 3H_1 + 2H_2 \ ,\quad
 H_5 = H_1+H_2 \ , \quad  H_6 = 2H_1+H_2 \ee
(compare with \eqref{dual}).

Note that for \,$H_k$\, also holds:
\eqn{cons} \l (H_k) \,=\, (\l,\a^\vee_k) \ , \qquad \forall\,\l\in\ch^* \ , \quad k=1,\ldots,6.\ee

\subsection{Structure theory of the real form}

The split real form of ~$G_2$~ is denoted as ~$G'_2\,$, sometimes as ~$G_{2(2)}\,$.
This real form has quaternionic discrete series representations.  We can use
the same basis (but over $\bbr$) and the same root system.

The Iwasawa decomposition of the real split form ~$\cg \equiv G'_2\,$,~ is:
\eqn{iwa32} \cg ~=~ \ck \oplus \ca_0 \oplus \cn_0 \ , \ee
the Cartan decomposition is:
\eqn{car32} \cg ~=~ \ck \oplus \cq , \ee
where we use: the maximal compact subgroup ~$\ck \cong
su(2) \oplus su(2)$, ~$\dim_\bbr\,\cq = 8$,
 ~$\dim_\bbr\,\ca_0 = 2$, ~$\cn_0 =  \cn^+_0\,$, or
 ~$\cn_0 =  \cn^-_0 \cong \cn^+_0\,$, ~$\dim_\bbr\,\cn_0^\pm =6$.

   Since ~$\cg$~ is maximally split, then
the centralizer ~$\cm_0$~ of ~$\ca_0$~ in ~$\ck$~ is zero, thus, the minimal parabolic ~$\cp_0$~
and the corresponding Bruhat decomposition are:
\eqn{min32} \cp_0 ~=~ \ca_0 \oplus \cn_0 \ , \qquad \cg ~=~ \ca_0 \oplus \cn^+_0\oplus \cn^-_0 \ee

The importance of the parabolic subgroups comes from the fact that
the representations induced from them generate all (admissible)
irreducible representations of the group under consideration \cite{Lan,Zhea,KnZu}.

\section{Verma modules and singular vectors}

Let us recall that a ~{\it Verma module} ~$V^\L$~ is defined as
the highest weight module over ~$\cgc$~ with highest weight ~$\L \in \ch^*$~ and
highest weight vector ~$v_0 \in V^\L$, induced from the
one-dimensional representation ~$V_0 \cong \bbc v_0$~ of
~$U(\cb)$~, where ~$\cb  = \ch \oplus \cg_+$~ is a Borel
subalgebra of ~$\cgc$, such that:
\eqnn{indb}
& &X ~v_0 ~~=~~ 0 , \quad \forall\, X\in \cg_+ \cr
&&H ~v_0 ~~=~~ \L(H)~v_0\,, \quad \forall\, H \in \ch \eea

Verma modules are generically irreducible. A Verma  module ~$V^\L$~ is
reducible \cite{BGG} iff there exists a root ~$\b \in\D^+$~ and ~$m\in\bbn$~
such that
\eqn{red} (\L + \r~, ~\b^\vee) ~=~  m \ee
holds, where ~$\r = {1 \over 2}\sum_{\a \in \D^+}~\a$~.
If \eqref{red} holds then the reducible Verma module $V^\L$  contains an invariant submodule
which is also a Verma module ~$V^{\L'}$~ with shifted weight ~$\L'=\L-m\b$.
This statement is equivalent to the fact that $V^\L$ contains a
~{\it singular vector}~
~$v_s \in V^\L$, such that ~$v_s ~\neq ~\xi v_0\,$, ($0\neq\xi\in\bbc$),
and~:
\eqnn{sing}
&& X ~v_s ~~=~~ 0 , \quad \forall\, X\in \cg_+ \cr
&&H ~v_s ~~=~~ \L'(H) ~v_s\,, \quad
\L' ~=~ \L - m\b, ~~\forall\,
H \in \ch \eea

%We now consider our main tool - the singular vectors of
%Verma modules ~$V^\L$~ of highest weight $\L$.
   The general reducibility conditions \eqref{red} for $V^\L$
   spelled out for the six positive roots in our
situation are: \eqna{reda}
m_1 \,&=&\, m_1(\L) \,\doteq\,  \L(H_1) + 1  \,\in\bbn \ ,  \\
m_2 \,&=&\, m_2(\L) \,\doteq\, \L(H_2) + 1   \,\in\bbn \ ,   \\
m_3 \,&=&\, m_3(\L) \,\doteq\, \L(H_3) + 4 = 3m_1 + m_2   \,\in\bbn \ ,  \\
m_4 \,&=&\, m_4(\L) \,\doteq\,  \L(H_4) + 5 = 3m_1 + 2m_2  \,\in\bbn \ ,\\
m_5 \,&=&\, m_3(\L) \,\doteq\, \L(H_5) + 2 = m_1 + m_2   \,\in\bbn \ ,  \\
m_6 \,&=&\, m_4(\L) \,\doteq\,  \L(H_6) + 3 = 2m_1 + m_2  \,\in\bbn \ .
\eena
The singular vectors corresponding to these cases   are:
\eqna{svsv}
 v^{\a_1,m_1} \,&=&\,  (X^-_1)^{m_1}\, v_0 \,,  \quad m_1\in\bbn\ ,\\
v^{\a_2,m_2} \,&=&\,  (X^-_2)^{m_2}\, v_0 \,, \quad m_2\in\bbn\ , \\
v^{\a_3,m_3} \,&=&\,  \sum_{k=0}^{m_3}\ g_{3k}\
(X_1^-)^{m_3-k}\, (X_2^-)^{m_3} \, (X_1^-)^{k}\, v_0 \,,
\quad m_3\in\bbn\ , \\
g_{3k} \,&=&\, \bca  g_3\, (-1)^k\, {m_3 \choose k} \, { m_1 \over
m_1 -k}\ , \,\,\,&\,\, m_1 \notin \{\, 0,\ldots, m_3\,\}  \cr
g_3 \,\d_{k,m_1} \ , \,\,\,&\,\, m_1 \in \{\, 0,\ldots, m_3\,\}
\eca \nn\\
v^{\a_5,m_5} \,&=&\,  \sum_{k=0}^{3m_5} \ g_{5k}
\ (X^-_2)^{3m_5-k} \ (X^-_1)^{m_5} \ (X^-_2)^k \  v_0 \ , \quad m_5\in\bbn\\\
g_{5k} \ &=&\  \bca  (-1)^k \ g_5 \ \left({3m_5 \atop k}\right) \  { m_2 \over
m_2 -k}\ , \,\,\,&\,\, m_2 \notin \{\, 0,\ldots, m_5\,\}  \cr
g_5\, \d_{k,m_2} \ , \,\,\,&\,\, m_2 \in \{\, 0,\ldots, m_5\,\}
\eca \nn \\
v^{\a_4,m_4} &=&  \sum_{k=0}^{m_4}  \sum_{i=0}^{2m_4} \sum_{j=0}^{2m_4-i} g_{4ijk}
(X_2^-)^{j} (X_1^-)^{m_4-k} (X_2^-)^{2m_4-i-j}  (X_1^-)^{k} (X_2^-)^{i}v_0 \nn\\
\\ %\quad m_4\in\bbn\
v^{\a_6,m_6} &=&  \sum_{k=0}^{3m_6}  \sum_{i=0}^{2m_6} \sum_{j=0}^{2m_6-i} g_{6ijk}
(X_1^-)^{j} (X_2^-)^{3m_6-k} (X_1^-)^{2m_6-i-j}  (X_2^-)^{k} (X_1^-)^{i}v_0 \nn\\
%\quad m_4\in\bbn\
 \eena
(Note that in each of the six cases \eqref{svsv} only the relevant $m_j$ must
be a natural number (as displayed).)
Formulae \rf{svsv}{a,b} are general for any simple root \cite{Dix},\cite{Dob},
while \rf{svsv}{c,d} were given first in \cite{Dos}.

Certainly, \eqref{red} may be fulfilled for several positive roots, even
for all of them if \eqref{red} is fulfilled for the two simple roots.

\section{Classification of  $G_2$ Verma modules}
% and P$_0$-induced ERs}
 \label{class}

 Here we classify the Verma modules over \,$\cg^\bbc=G_2)$.
 This also provides the classification of
 the $P_0$-induced ERs since the restricted Weyl group \ $W(\cg,\ca_0)$\  related to the
 minimal parabolic subalgebra \ $\cp_0=\ca_0\cn_0\,$, cf. \eqref{min32}, is isomorphic to
 the Weyl group \ $W(\cg^\bbc,\ca_0^\bbc)$ (since $\cg$ is maximally split).

The classification is done as follows.   We group the
reducible Verma modules  related by nontrivial embeddings
%with the same Casimirs
in sets called ~{\it multiplets} \cite{Dobmul,Dob}. These multiplets
% corresponding to fixed values of the Casimirs
may be depicted as a connected graph, the
vertices of which correspond to the reducible Verma modules and the lines
between the vertices correspond to the embeddings. %intertwining operators.
The explicit parametrization of the multiplets and of their Verma modules is
important for understanding of the situation.

The classification can be summarized as follows. There are four main types of
multiplets of reducible Verma modules:
\begin{itemize}
\item{}\,type \,$\cf_{m_1,m_2}\,$, ($m_1,m_2\in\bbn$);
\item{}  type ~$\cn$~ with six subtypes:
~$\cn_k$, $k=1,\ldots,6$;
\item{}\,type \,$\cl_{m}\,$, with two subtypes $\cl^k_{m}\,$, $k=1,2$, ($m\in\bbn$);
\item{}\,type \,$\cs_{j,k}\,$, ~$m_j,m_k\in\bbn$, \ $1\leq j<k \leq 6$, ~$(j,k)\neq (1,2)$;
%\item{}\,type \,$\car$\, with four subtypes: \,$\car_k$, $k=1,2,3,4$;
\end{itemize}

\md

Multiplets of type \,$\cf_{m_1,m_2}$\, are parametrized by two
natural numbers \,$m_1,m_2\,$. They are given in Fig. 1A and Fig. 1B
where we have given the multiplets
in two ways: in Fig. 1A the Verma modules are given by their
highest weights, while on Fig. 1B  they are given by the two Dynkin labels. %Harish-Chandra parameters.
 In Fig. 1A we have indicated   w.r.t. which reflection is the embedding,
  on Fig. 1B  we have given the weight $m\b$ of the embedding. All Verma modules of
these multiplets, except $V^{\L'}$, are reducible.
% and their %weights are given explicitly as follows:
%\eqnn{fdir}  &&
%\L_{m_1,m_2} \ , \quad m_k = m_k(\L_{m_1,m_2}) \in \bbn, \,\,\forall
%k \ , \\ && \L_{m_1,m_2}^i = \L_{m_1,m_2} - m_i\a_i\ , \,\,\,i=1,2 \
%,\nn\\ && \L_{m_1,m_2}^{ij} = \L_{m_1,m_2} - m_i\a_i - m_{j+2}\a_j \
%, \,\,\,(i,j)=(1,2),(2,1)\ , \nn\\ && \L_{m_1,m_2}^{iji} = \L_{m_1,m_2}
%- (m_i+m_{i+2})\a_i - m_{j+2}\a_j \ , \,\,\,(i,j)=(1,2),(2,1)\  \nn\eea
%The weights of the irreducible modules $V^{\L^+}$ are: \,$\L^+_{m_1,m_2} \ =\
%\L_{m_1,m_2} - (m_1+m_3)\a_1 - (m_2+m_4)\a_2 = \L_{m_1,m_2} - 2m_4\a_1 -
%m_3\a_2\,$.
Note that only embeddings which are not compositions
of other embeddings are given on the Figures.

\fig{}{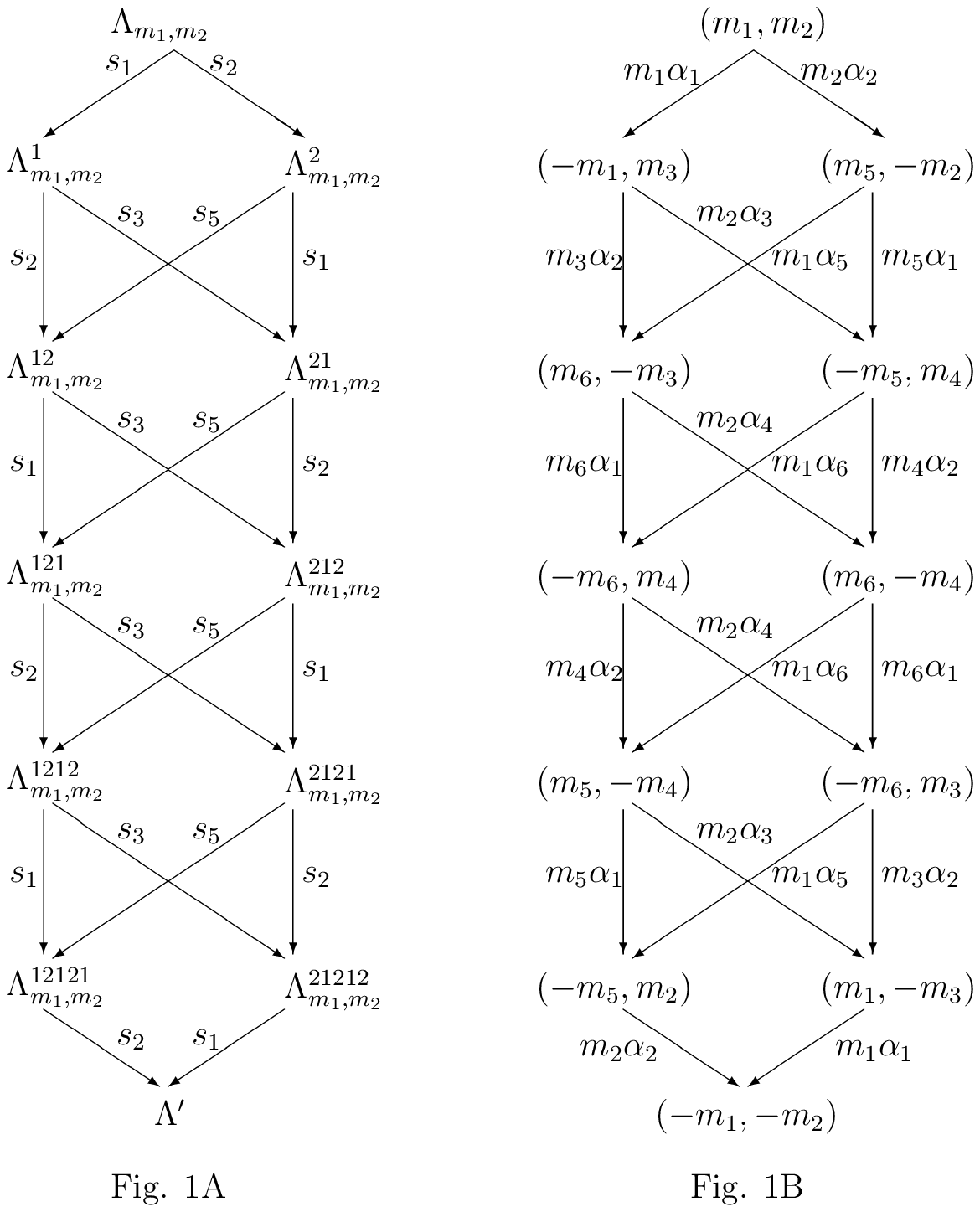}{8cm}
%\fig{}{vectorg2b.eps}{6cm}
%\np

We note also some additional relations using notation from Fig. 1:
\eqnn{weyll}
&& \L^1 = \L-m_1\a_1, \quad  \L^2 = \L-m_2\a_2, \quad
\L^{121} = \L-m_3\a_3, \quad  \L^{212} = \L-m_5\a_5, \nn\\
&& \L^{12121} = \L-m_6\a_6, \quad  \L^{21212} = \L-m_4\a_4 \eea

\md

Multiplets of type ~$\cn$~  are given as follows. Fix $k=1,\ldots,6$
to fix a subtype ~$\cn_k\,$.~ Then the multiplets of this subtype
are parametrized by the natural number $m_k$ and are given as follows:
\eqn{typa}
 V^{\L_k} ~\longrightarrow~   V^{\L_k - m_k\a_k} \ ,
\quad  m_k(\L_k)= m_k\in\bbn \ , \quad m_j(\L_k)\notin\bbn, ~~j\neq k\ . \ee %&\typa a\cr
%&& L_{\L_k} ~=~ V^{\L_k} / V^{\L_k - m_k\a_k} \ .\eena
Note that we are using the convention that the arrows point to the embedded modules.
The modules $V^{\L_k - m_k\a_k}$ are irreducible.

\md

For the multiplets of type \,$\cl$\,
there are two subtypes \,$\cl_k\,$,  \ $k=1,2$\,
each parametrized by a natural number $m$ and  given as follows.

The multiplets of subtype \,$\cl_1\,$\, are given on  Fig. 2 below.
In each multiplet there are six Verma modules which we give by the Dynkin labels.
We also give the weights of the singular vectors between the Verma modules.
The last module on the right is irreducible.

\bigskip

\fig{}{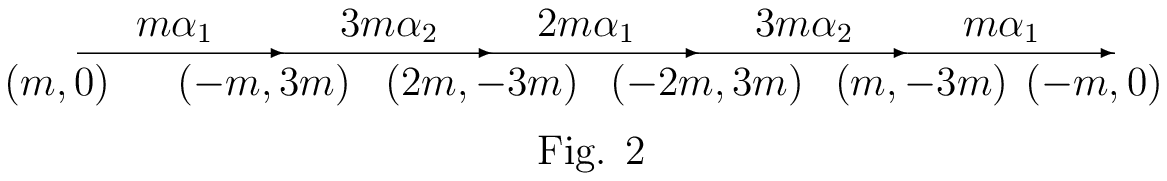}{8cm}

The multiplets of subtype \,$\cl_2\,$\, are given on Fig. 3 below.
They are similar to subtype \,$\cl_1\,$, e.g.,
also here the last module on the right is irreducible.

\bigskip

 \fig{}{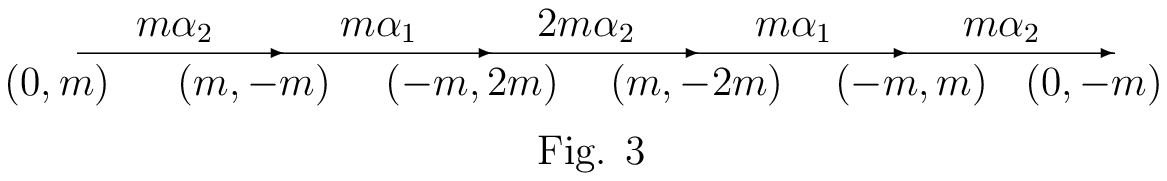}{8cm}\md

\md

For the lack of space we give only some examples of the multiplets of type $\cs_{j,k}$.

Multiplets of type \,$\cs_{1,4}$\, are parametrized by two natural numbers \,$m_1,m_4$\,
 so that \,$m_2 = \ha (m_4-3m_1) \notin\bbz$; then also $m_3,m_5,m_6\notin\bbz$.
 They are given in the Fig. 4 below,
where as above we have given the multiplet in two ways,
and again the parametrizing numbers $m_1,m_4$ are related to the Verma module $V^{\L^s}$
on the top: \,$m_k = m_k(\L^s)$, $k=1,4$.
The Verma modules of these multiplets, except $V^{\L''}$, are reducible
and their weights are given explicitly as follows:
\eqn{fdis}
\L^s \ , \quad m_k = m_k(\L^s) \in \bbn, \,\,k=1,4,
 \quad \L^s_k = \L^s - m_k\a_k\ , \,\,\,k=1,4 \ .
  \ee
The weights of the irreducible modules $V^{\L'}$ are:
\,$\L' = \L^s - m_1\a_1 - m_4 \a_4\,$. \md
%$$\includegraphics[bb=-10 0 500 100,
% width=14cm]{vector13}$$

\md

\fig{}{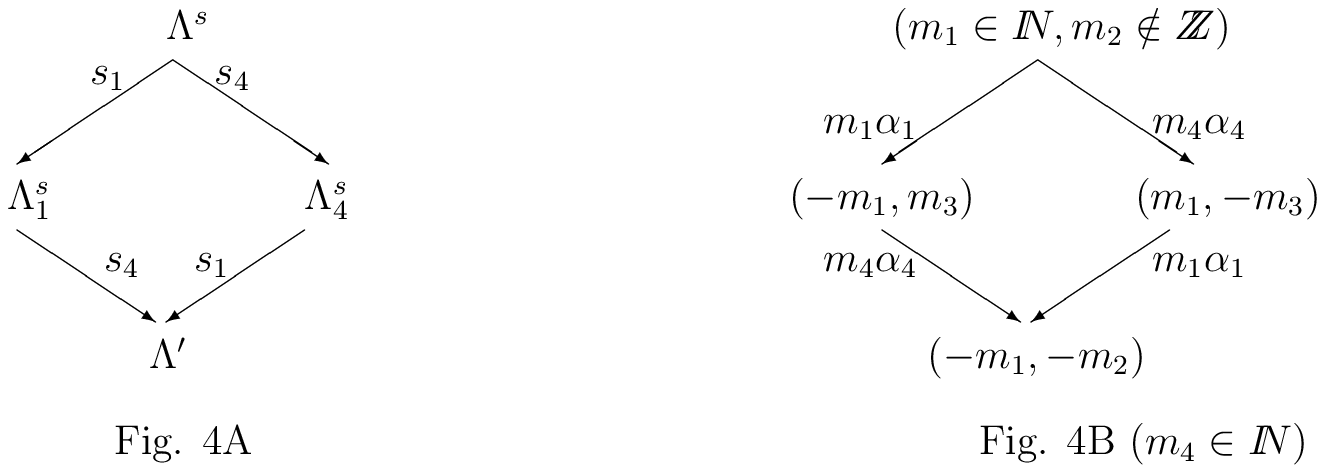}{11cm}

\md
Multiplets of type \,$\cs_{2,6}$\, are parametrized by two natural numbers \,$m_2,m_6$, so that
~$m_1 = \ha (m_6-m_2) \notin\bbz$.
Here we have two subcases: a) ~$3m_1 = \trha (m_6-m_2) \notin\bbz$, then also $m_3,m_4,m_5\notin\bbz$;
b) ~$3m_1 = \trha (m_6-m_2) \in\bbn$, then $m_3,m_4\in\bbn$, $m_5\notin\bbz$.\\
 Subcase a) is given in the Fig. 6 below, and we omit comments since this case is similar to Fig. 5.
\fig{}{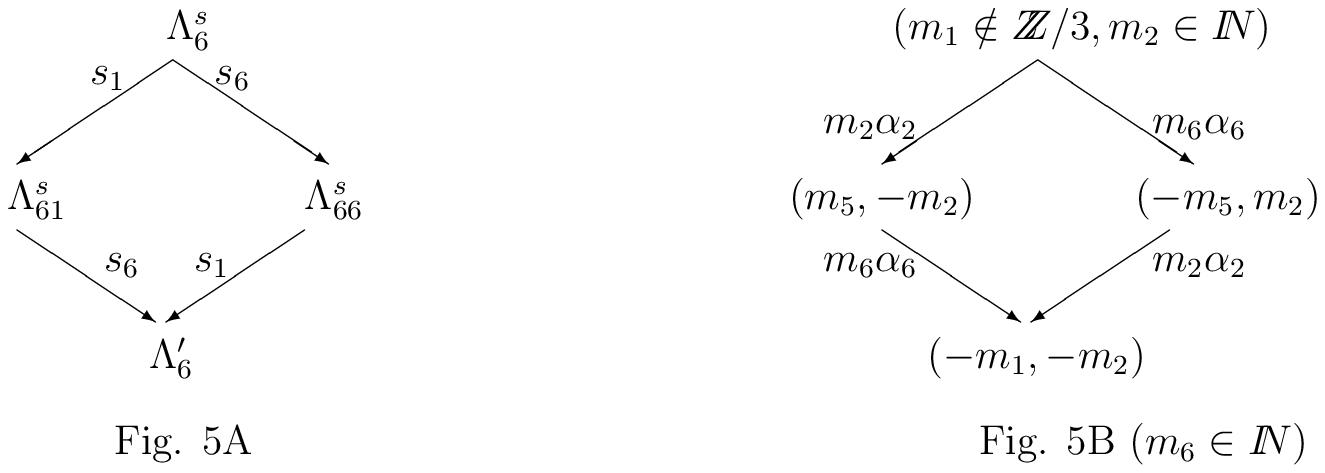}{11cm}

Subcase b) contains eight Verma modules as given in Fig. 6 below. Here we use a mixture of notation since
the Dynkin labels may be seen from comparing Fig. 1 and Fig. 2. Here there are two irreducible modules:
~ $\L'$~ and ~$\Lambda^{21212}$ (the second one since $m_1\notin\bbn$).\md
\fig{}{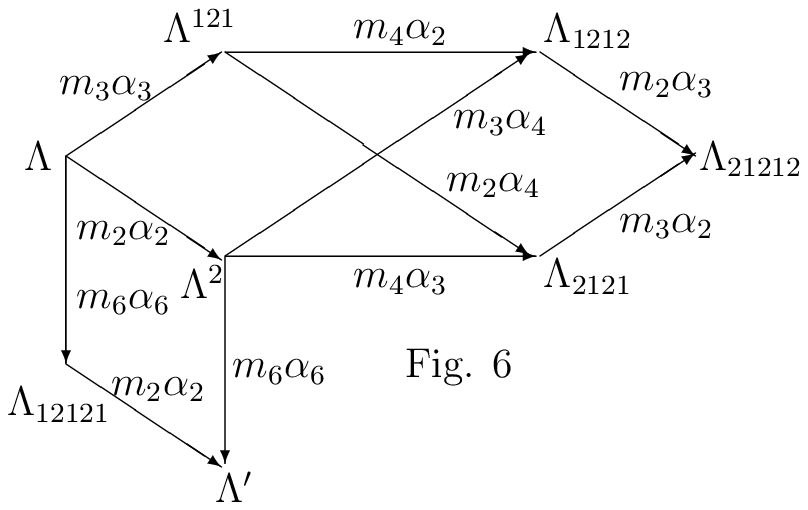}{7cm}
%\fig{}{vectorg7.eps}{7cm}

\md
Multiplets of type \,$\cs_{2,4}$\, are parametrized by two natural numbers \,$m_2,m_4$, so that
~$m_1 = \third (m_4-2m_2) \notin\bbn$, while $m_3 = m_4-m_2\in\bbn$, ~~$m_5,m_6\notin\bbn$.
These multiplets are given in the Fig. 7 below: \md
\fig{}{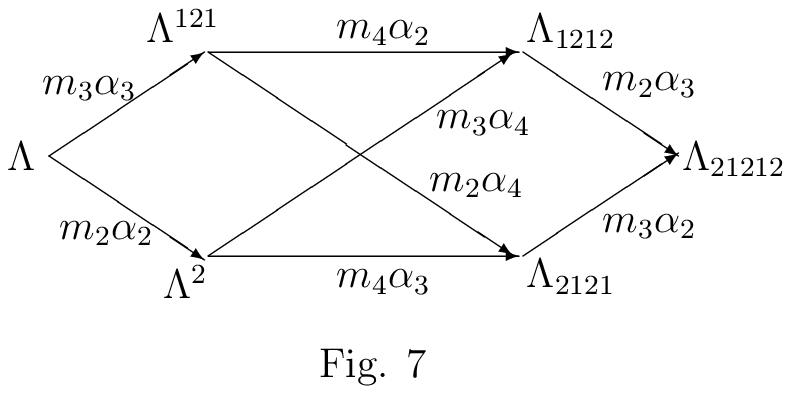}{7cm}
\md
\nt Note that this multiplet is the standard $sl(3)$ multiplet with $\a_2,\a_3$ playing the role
of the two simple roots. Incidentally, it is a submultiplet of the previous case.

\bigskip

\section*{Acknowledgments.}\nt
 The author would like to thank the Organizers for the kind
hospitality and invitation to present a talk at the XXXII
International Colloquium on Group Theoretical Methods in Physics
(Prague, July 2018). The author has received partial support from
Bulgarian NSF Grant DN-18/1, from COST Action MP1405 and from PHC
Rila.

\end{document}